\documentclass[12pt, reqno]{amsart}
\usepackage{graphicx}

\newtheorem{theorem}{Theorem}[section]
\newtheorem{proposition}[theorem]{Proposition}
\newtheorem{lemma}[theorem]{Lemma}
\newtheorem{corollary}[theorem]{Corollary}

\theoremstyle{definition} 
\newtheorem{definition}[theorem]{Definition}

\newtheorem{remark}[theorem]{Remark}
\newtheorem{question}[theorem]{Question}

\newcommand{\C}{\mathbb C} 
 \newcommand{\T}{\mathbb T}

\DeclareMathOperator{\sys}{{\rm sys}}

\DeclareMathOperator{\pisys}{\sys\!\pi}

\DeclareMathOperator{\vol}{{\rm vol}}

\DeclareMathOperator{\area}{{\rm area}}
\DeclareMathOperator{\length}{{\rm length}}

\DeclareMathOperator{\gmetric}{{\mathcal G}}

\def\SR{{\rm SR}}

\def\tens{{\rm tens}}

\def\LL{{L_L}}

 \def\gen {{\it g}}

\def\ie {{\it i.e.\ }} 
\def\eg {{\it e.g.\ }} 
\def\cf {\hbox{\it cf.\ }}

\numberwithin{equation}{section} 
\numberwithin{figure}{section}

\begin{document}

\author[M.~Katz]{Mikhail G. Katz$^{*}$} \address{Department of
Mathematics and Statistics, Bar Ilan University, Ramat Gan 52900
Israel} \email{katzmik@math.biu.ac.il} \thanks{$^{*}$Supported by the
Israel Science Foundation (grants no.\ 620/00-10.0 and 84/03)}

\author[S.~Sabourau]{St\'ephane Sabourau} \address{Laboratoire de
Math\'ematiques et Physique Th\'eorique, Universit\'e de Tours, Parc
de Grandmont, 37400 Tours, France}

\address{Mathematics and Computer Science Department, Saint-Joseph's
University, 5600 City Avenue, Philadelphia, PA 19131, USA}

\email{sabourau@gargan.math.univ-tours.fr}

\title{Hyperelliptic surfaces are Loewner$^1$}


\begin{abstract}
We prove that C. Loewner's inequality for the torus is satisfied by
conformal metrics on hyperelliptic surfaces $X$, as well.  In genus 2,
we first construct the Loewner loops on the (mildly singular)
companion tori, locally isometric to $X$ away from Weierstrass points.
The loops are then transplanted to $X$, and surgered to obtain a
Loewner loop on~$X$.  In higher genus, we exploit M.~Gromov's area
estimates for $\varepsilon$-regular metrics on $X$.
\end{abstract}

\footnotetext[1]{\large {\em Proceedings Amer.~Math.~Soc.}, to appear}

\maketitle

\tableofcontents

\section{Hermite constant and Loewner surfaces}

The systole, $\pisys_1(\gmetric)$, of a compact non simply connected
Riemannian manifold~$(X,\gmetric)$ is the least length of a
noncontractible loop $\gamma\subset X$:
\begin{equation}
\label{ps1}
\pisys_1(\gmetric)=\min_{[\gamma]\not= 0\in \pi_1(X)}\length(\gamma).
\end{equation}
This notion of systole is apparently unrelated to the systolic arrays
of~\cite{Ku}.  We will be concerned with comparing this Riemannian
invariant to the total area of the metric, as in the Loewner
inequality~\eqref{(1.1)} for the torus.

The Hermite constant, denoted $\gamma_n$, can be defined as the
optimal constant in the inequality
\begin{equation}
\label{11}
\pisys_1(\T^n)^2 \leq \gamma_n \vol(\T^n)^{2/n},
\end{equation}
over the class of all {\em flat\/} tori $\T^n$.  Here $\gamma_n$ is
asymptotically linear in~$n$, \cf \cite[pp.~334, 337]{LLS}.  The
precise value is known for small $n$, \eg one has~$\gamma_2=
\frac{2}{\sqrt{3}}$, $\gamma_3 = 2^{\frac{1}{3}}$.  An inequality of
type~\eqref{11} remains valid in the class of {\em all\/} metrics, and
in fact on a more general manifold, but with a nonsharp constant on
the order of $n^{4n}$ \cite{Gr1}.

Around 1949, Charles Loewner proved the first systolic inequality, \cf
\cite{Pu}.  He showed that every Riemannian metric $\gmetric$ on the
torus~$\T^2$ satisfies the inequality
\begin{equation}
\label{(1.1)}
\pisys_1(\gmetric)^2\le \gamma_2 \area(\gmetric),
\end{equation} 
while a metric satisfying the boundary case of equality in
\eqref{(1.1)} is necessarily flat, and is homothetic to the quotient
of $\C$ by the lattice spanned by the cube roots of unity.

Interesting extensions of the Loewner inequality are studied in
\cite{Am}.  Higher dimensional optimal generalisations of the Loewner
inequality are studied in \cite{BanK, Ka3, BaKa, IK, BCIK2, Sa2}.  The
defining text for this material is \cite{Gr3}, with more details in
\cite{Gr1,Gr2}, complemented by \cite{bab}.  See also the recent
survey \cite{CK}, as well as \cite{KL, KR}.

\begin{definition} 
We say that a metric~$\gmetric$ on a surface is {\em Loewner\/} if it
satisfies the Loewner inequality~\eqref{(1.1)}.
\end{definition}

\begin{question}
\label{21}
It follows from Gromov's estimate \eqref{1.1} that every metric on an
orientable surface $\Sigma_\gen$ of genus $\gen$ is Loewner if $\gen >
50$.  This result is improved in \cite{KS2} to $\gen \geq 20$.  Can
the genus assumption be removed altogether?
\end{question}

Note that a similar question for Pu's inequality \cite{Pu} has an
affirmative answer.  The generalisation is immediate from Gromov's
inequality~\eqref{tqi}.  Namely, every surface $X$ which is not a
2-sphere satisfies the inequality
\begin{equation}
\label{pu3}
\pisys_1(X)^2 \leq \frac{\pi}{2} \area(X),
\end{equation}
where the boundary case of equality in \eqref{pu3} is attained
precisely when, on the one hand, the surface $X$ is a real projective
plane, and on the other, the metric is of constant Gaussian curvature.

We will prove the following result toward answering Question \ref{21},
\cf Theorem~\ref{31}.

\begin{theorem} 
Every metric on an orientable surface is Loewner if it lies in a
hyperelliptic conformal class. In particular, every metric on the
genus~2 surface is Loewner.
\end{theorem}

While the precise value of the systolic ratio in genus 2 is unknown in
the class of all metrics, an optimal systolic inequality for CAT(0)
metrics does exist \cite{KS3}.  The relevant literature and basic
estimates are reviewed in Section~\ref{one}.  We will state the main
theorem in more detail in Section~\ref{three} and prove it for genus 3
or more.  We will complete the proof in the genus 2 case in
Section~\ref{four}.

\section{Basic estimates}
\label{one}

\begin{definition}
The {\em systolic ratio} of a metric~$\gmetric$ on a closed $n$-manifold is
defined as
\[
\SR(\gmetric) = \frac{\pisys_1(\gmetric)^n}{\vol_n(\gmetric)}.
\]
The {\em conformal systolic ratio}, denoted $\SR_c(X)$, of a closed
$n$-manifold~$X$ with a chosen conformal class of metrics, is defined as
\[
\SR_c(X) = \sup_{\gmetric} \left\{ \SR(\gmetric) : \gmetric \mbox{ a
conformal metric on } X \right\}.
\]
The supremum of the conformal systolic ratio over all the conformal
structures of~$X$ is called the {\em optimal systolic ratio}.  It is
denoted by~$\SR(X)$.
\end{definition}
M.~Gromov~\cite[p.~50]{Gr1} (\cf \cite[Theorem 4, part (1)]{Ko})
proved a general estimate which implies that if $\Sigma_\gen$ is a
compact orientable surface of genus $\gen$ with a Riemannian metric,
then
\begin{equation}
\label{1.1}
\SR(\Sigma_\gen) < \frac {64} {4\sqrt{\gen} + 27} .
\end{equation}
Thus the optimal systolic ratio tends to $0$ as the genus increases
without bound.

\begin{remark}
\label{21b}
It was shown in \cite{Gr1} (see also~\cite{KS2} and~\cite{bab}) that
asymptotically the optimal systolic ratio of a surface of genus $\gen$
behaves as~$C \frac{(\log \gen)^2}{\gen}$, as $\gen \to \infty$.
\end{remark}

Another helpful estimate is found in \cite[Corollary~5.2.B]{Gr1}.
Namely, every aspherical compact surface $(\Sigma,\gmetric)$ admits a metric
ball $B=B_p\left(\tfrac{1}{2}\pisys_1(\gmetric)\right) \subset \Sigma$ of
radius $\tfrac{1}{2}\pisys_1(\gmetric)$, which satisfies
\begin{equation}
\label{tqi}
\pisys_1(\gmetric)^2 \leq \frac{4}{3}\area(B).
\end{equation}
Furthermore, whenever a point $x\in \Sigma$ lies on a two-sided loop
which is minimizing in its free homotopy class, the metric ball $
B_x(r)\subset \Sigma $ of radius $r\leq \frac{1}{2}\pisys_1(\gmetric)$
satisfies the estimate 
\begin{equation}
\label{24}
2r^2 < \area \left( B_x(r) \right).
\end{equation}

\section{Hyperelliptic surfaces and $\varepsilon$-regularity}
\label{three}

Recall that a Riemann surface $X$ is called hyperelliptic if it admits
a degree 2 meromorphic function, \cf \cite[p.~60-61]{Mi} as well as
\cite[Proposition 4.11, p.~92]{Mi}.  The associated ramified double
cover
\begin{equation}
\label{31b} 
Q: X\to S^2
\end{equation}
over the sphere $S^2$ is conformal away from the $2\gen+2$ ramification
points, where~$\gen$ is the genus of $X$.  Its deck transformation $J: X
\to X$ is called the {\em hyperelliptic involution\/}.  Such a
holomorphic involution, if it exists, is uniquely characterized by the
property of having precisely~$2\gen+2$ fixed points.  The fixed points of
$J$ are called Weierstrass points.  Their images under the map~$Q$ of
\eqref{31b} will be referred to as {\em ramification points}.

\begin{theorem} 
\label{31}
Let $(X,\gmetric)$ be an orientable surface, where the
metric~$\gmetric$ belongs to a hyperelliptic conformal class.  Then
$(X,\gmetric)$ is Loewner.
\end{theorem}

Since every genus 2 surface is hyperelliptic \cite
[Proposition~III.7.2, page~100] {FK}, we obtain the following
corollary.

\begin{corollary}
\label{65}
Every metric on the genus 2 surface is Loewner.
\end{corollary}

Note that this is the first improvement, known to the authors, on
Gromov's 3/4 bound \eqref{tqi} in over 20 years, for surfaces of genus
below~50, \cf Question~\ref{21}.  No extremal metric has as yet been
conjectured in this genus, but it cannot be flat with conical
singularities~\cite{Sa}.  The best available lower bound for the
optimal systolic ratio in genus 2 can be found in
\cite[section~2.2]{CK}.  

For genus $\gen \geq 3$, our Theorem~\ref{31} follows from
Proposition~\ref{s>3}, \cf Remark~\ref{21b} and~\cite{Kon}.  Let us go
over the definitions of~\cite[5.1]{Gr1}.

\begin{definition}
Given a Riemannian metric~$\gmetric$ on $X$, the {\em tension\/} of a
noncontractible loop~$\gamma$ based at~$x\in X$ is the upper bound of
all~$\delta>0$ such that there exists a free homotopy of~$\gamma$
which diminishes the length of~$\gamma$ by~$\delta$.  The tension is
denoted~$\tens_{\gmetric}(\gamma)$.
\end{definition}

\begin{definition}
The height~$h_{\gmetric}(x)$ of $x\in X$ is defined as the lower bound
of the tensions of the noncontractible loops based at~$x$. Note that
the height function is $2$-Lipschitz.

A metric~$\gmetric$ is said to be $\varepsilon$-regular with
$\varepsilon < \pisys_1(\gmetric)$ if its height
function~$h_{\gmetric}$ is bounded from above by~$\varepsilon$.
\end{definition}

\begin{lemma}
\label{reg}
Given any conformal metric $\gmetric$ on~$(\Sigma_\gen,J)$ and a real
number $\delta>0$, there exists an~$\varepsilon$-regular,
$J$-invariant, conformal metric~$\bar{\gmetric}$ with a systolic ratio
at least~$\SR(\gmetric) - \delta$.
\end{lemma}

\begin{proof}
The proof appears in~\cite[$5.6.C''$]{Gr1} in the general case, and
proceeds by a (finite) sequence of modifications of the metric in
suitable small disks, while staying in the same conformal class.  Note
that averaging the metric by $J$ improves the systolic ratio and does
not change the conformal class.  Thus we can assume that $\gmetric$ is
$J$-invariant.  To adapt the proof to our situation, we perform the
modifications in a $J$-invariant way.
\end{proof}

\begin{proposition}
\label{s>3}
Every hyperelliptic surface $(\Sigma_s, J)$ of genus~$\gen$ satisfies
the estimate
\[
\SR_c(\Sigma_\gen) \leq \frac{4}{\gen+1}.
\]
\end{proposition}

\begin{proof}
Let $\gmetric$ be a conformal metric on the surface.  By
Lemma~\ref{reg}, we can assume that the metric~$\gmetric$ is
$J$-invariant and~$\varepsilon$-regular.  Since the metric is
$J$-invariant, the distance between any pair of the $2\gen+2$ Weierstrass
points is at least $\frac{1}{2} \pisys_1(\gmetric)$.  Thus, the disks
of radius $R=\frac{1}{4} \pisys_1(\gmetric)$ centered at the
Weierstrass points are disjoint.  From~\cite[5.1.B]{Gr1} (see
also~\cite{He}), the area of such a disk on an~$\varepsilon$-regular
surface is at least
\[
2R^2 - \varepsilon =\frac{1}{8} \pisys_1(\gmetric)^2 - \varepsilon,
\]
\cf \eqref{24}.
Therefore, we have
\[
\area(\gmetric) \geq (2 \gen +2) \left(\tfrac{1}{8}
\pisys_1(\gmetric)^2 - \varepsilon \right).
\]
Since this is true for all $\varepsilon>0$, the proposition is proved.
\end{proof}

\section{Proof of Theorem~\ref{31} in genus two}
\label{four}

Let $X$ be a genus~$2$ surface. Recall that $X$ has a hyperelliptic
involution~$J$ with $6$ Weierstrass points.

The idea of the proof of Theorem~\ref{31} in genus~$2$ is to apply the
Loewner inequality to certain {\em companion tori\/} of $X$, and to
surger the resulting loops so as to obtain a Loewner loop on $X$.  We
may need the following lemma.

\begin{lemma} 
\label{lem:simple}
Let $\T^2$ be a torus endowed with a metric invariant under its
hyperelliptic involution $J_{T^2}$, with conical singularities with
total angle less than~$2 \pi$ around each.  Then the image of a
systolic loop of~$\T^2$ in $S^2$ under the hyperelliptic projection is
a simple loop.
\end{lemma}

\begin{proof}
Let $\gamma \subset \T^2$ be a systolic loop.  Since $J_{T^2}$ induces
minus the identity homomorphism on~$\pi_1(\T^2)$, the loops $\gamma$
and $-J_{T^2}(\gamma)$ are homotopic.  Under the hypotheses of our
lemma, two homotopic systolic loops are necessarily disjoint.  Hence
the image of $\gamma$ on $S^2$ is simple.
\end{proof}

\begin{definition}
A {\em companion torus\/} $\T(a,b,c,d)$ of $X$ is a torus whose
ramification locus $\{a,b,c,d\}\subset S^2$ is a subset of the
ramification locus of $X$.
\end{definition}

As in the proof of Proposition~\ref{s>3}, we can assume that the
metric on $X$ is invariant under $J$ (see~\cite{BCIK1}).  Therefore
$\gmetric$ descends to a metric $\gmetric_0$, of half the area, on
$S^2$.  Let's choose four of the $6$ ramification points, say
$a,b,c,d\in S^2$.  Choose a double cover with ramification
locus~$\{a,b,c,d\}$, denoted
\[
\T^2(a,b,c,d)\to S^2.
\]
Pulling back the metric~$\gmetric_0$ to the torus $\T^2(a,b,c,d)$, we
obtain a metric of the same area as the surface $X$ itself.  This
metric on the torus is smooth away from the two remaining points,
where it has a conical singularity with total angle~$\pi$ around each.
Consider a Loewner loop
\[
\LL \subset \T^2(a,b,c,d)
\]
on this torus, \eg a systolic loop realizing \eqref{(1.1)}.  Let~$L$
be the projection of~$\LL$ to $S^2$.  The simple loop~$L\subset S^2$
separates the four points $a,b,c,d$ into two pairs, say $a,b$ on one
side and $c,d$, on the other.  If the lift of~$L$ to~$X$ closes up, we
obtain a Loewner loop on~$X$ and the theorem is proved.  Thus, we may
assume that the following three equivalent conditions are satisfied:

\begin{enumerate}
\item
the lift of $L$ to $X$ does not close up;
\item
the inverse image $Q^{-1}\left( L \right) \subset X$ under $Q$
of~\eqref{31b} is connected;
\item
the loop $L$ surrounds precisely $3$ ramification points of~$Q$.
\end{enumerate}

The last condition is equivalent to the first two since every based
loop is homotopic, in the complement of the ramification points, to a
composition of some loops from a finite collection of ``standard''
simple based loops, circling each of the ramification points.
Meanwhile, going once around such a standard loop clearly switches the
two sheets of the cover.

\begin{definition}
\label{33}
The simple loop $L$ partitions the sphere into two hemispheres,
$H_{+}$ and $H_{-}$, with $a,b,e\in H_+$ and $c,d,f\in H_-$ where
$a,b,c,d,e,f$ are the $6$ ramification points of~$Q$.
\end{definition}

Using a pair of companion tori, we will construct two loops on the
sphere, defining two distinct partitions of the ramification locus
into a pair of triples.  The basic example to think of is the case of
a centrally symmetric 6-tuple of points, corresponding for instance to
the curve
\[
y^2=x^5-x,
\]
and a pair of generic great circles, such that each of the four digons
contains at least one ramification point.  We now construct a
companion torus $\T(a,b,e,f)$.

Consider a Loewner loop $\LL'\subset \T^2(a,b,e,f)$, and its
projection $L' \subset S^2$.  If its lift to~$X$ closes up, the
theorem is proved.  Therefore assume that the lift of~$L'$ to~$X$ does
not close up, \ie $L'$ surrounds exactly $3$ ramification points.  Now
$L'$ separates the four points $a,b,e,f$ into two pairs.  Hence it
defines a different splitting of the six points into two triples.  The
connected components of $L' \cap H_+$ form a nonempty finite
collection of disjoint nonselfintersecting arcs $\alpha$.

Each arc $\alpha$ divides $H_+$ into a pair of regions homeomorphic to
disks.  Such regions are partially ordered by inclusion.  A minimal
element for the partial order is necessarily a digon.  Such a digon
must contain at least one ramification point of $Q$ (otherwise
exchange the two sides of the digon between the loops $L$ and $L'$, so
as to decrease the total number of intersections, or else argue as in
Lemma \ref{lem:simple}).  It is clear that there are at least two such
digons in $H_+$.  

Hence one of them, denoted $D\subset H_+$, must contain precisely one
of the 3 ramification points of $H_+$.  We now exchange the two sides
of~$D$ between the loops $L, L'$, obtaining two new loops $M, M'$.
Each of the new loops surrounds a nonzero even number of ramification
points.  Since 
\[
\length(M) + \length(M') = \length(L) + \length(L'),
\]
one of the loops~$M$ or $M'$ is no longer than Loewner.  Moreover, its
lift to $X$ closes up, producing a Loewner loop on $X$, as required.

\vfill\eject

\end{document}